\documentclass[12pt]{article}
\usepackage{amsmath, amsfonts}
\newcommand{\bee}{\begin{eqnarray*}}
\newcommand{\ene}{\end{eqnarray*}}
\newcommand{\beeq}{\begin{equation}}
\newcommand{\eneq}{\end{equation}}
\newtheorem{lem}{Lemma}[section]
\newcommand{\bel}{\begin{lem}}
\newcommand{\enl}{\end{lem}}
\newtheorem{defi}{Definition}[section]
\newcommand{\bef}{\begin{defi}}
\newcommand{\enf}{\end{defi}}
\newtheorem{exap}{Example}[section]
\newcommand{\beex}{\begin{exap}}
\newcommand{\enex}{\end{exap}}
\newtheorem{theo}{Theorem}[section]
\newcommand{\beth}{\begin{theo}}
\newcommand{\enth}{\end{theo}}
\newtheorem{prop}{Proposition}[section]
\newcommand{\bep}{\begin{prop}}
\newcommand{\enp}{\end{prop}}
\newtheorem{cor}{Corollary}[section]
\newcommand{\bec}{\begin{cor}}
\newcommand{\enc}{\end{cor}}
\newtheorem{rem}{Remark}[section]
\newcommand{\ber}{\begin{rem}}
\newcommand{\enr}{\end{rem}}

\usepackage{amsmath}
\usepackage{graphicx}
\usepackage{latexsym}
\begin{document}
\title{ 
  MLE'S BIAS PATHOLOGY,\\ 
 MODEL UPDATED MLE \\AND\\
WALLACE'S  MINIMUM MESSAGE LENGTH METHOD
}
\author{Yannis  G. Yatracos\\
Cyprus U. of Technology
}
\maketitle
\date{}

{\it e-mail:} yannis.yatracos@cut.ac.cy
 
\pagebreak 
 
 \begin{center}
{Summary}
\end{center}

\parbox{4.8in}

{
 The  inherent bias pathology of the maximum likelihood (ML)  estimation method is
confirmed 
 for
models  with  unknown parameters $\theta$ and $\psi$ when
MLE $\hat \psi$ is function of MLE $\hat \theta.$  To reduce $\hat \psi$'s bias
the likelihood equation to be solved for $\psi$   is updated using the
 model  for the data $Y$ in it.
Model updated (MU)  MLE, $\hat \psi_{MU},$ often  reduces either totally or partially $\hat \psi$'s  bias 
when estimating shape parameter $\psi.$  For the Pareto model
$\hat \psi_{MU}$  reduces
also $\hat \psi$'s variance.
The results  explain 
the  difference   that  puzzled R. A. Fisher, between biased $\hat \psi$ and the
unbiased estimate he obtained for two models with 
the ``2-stage procedure''.
MUMLE's implementation  is equivalent to the abandoned 2-stage procedure 
thus justifying its use.
 MUMLE  and Firth's  bias correcting likelihood are also obtained with the Minimum Message Length
method
thus motivating its use  in frequentist inference and, more generally,  model  updating with a prior distribution.


}
 
\bigskip
\vspace{.10in}
{\it Some key words:} \quad Bias, Likelihood 
equations,  Minimum Description Length Criterion, Minimum Message Length Method, 
Maximum likelihood, Model Updated MLE, Specification problem, Two-stage MLE

\date{ }


\pagebreak 

\section{Introduction}

\quad 
When data $x$ follows a model with density $f(x|\theta, \psi)$ and  parameters $\theta \in R^p ( p\ge 1),  \ \psi \in R, $   the maximum likelihood estimate (MLE)  $\hat \psi$ is  often biased when it depends on MLE $\hat \theta$ 
and the model is $\psi$-regular, i.e. the $\psi$-score's expectation  vanishes  for all $\theta, \psi.$ An alternative estimation  method for $\psi$ is thus motivated and  proposed.
The  model updated (MU)  maximum  likelihood principle (MLP) is used to obtain MUMLE,
$\hat \psi_{MU},$
 that reduces often $\hat \psi$'s bias and sometimes also its variance.
MUMLE and Firth's (1993) bias correcting likelihood are also obtained with the Minimum Message Length (MML) method (see, e.g. Wallace, 2005), i.e. by  either selecting  a $\psi$-prior to update $f(x|\hat \theta, \psi)$ 
and  obtain $\hat \psi_{MU}$ or decrease MLE's bias in general by updating $f(x|\theta, \psi)$ with a properly selected prior.

The results 
justify theoretically   Fisher's abandoned ``2-stage procedure'' that does not adhere to MLP and its implementation is equivalent to
MUMLE.
 When the MLE of a parameter has a distribution depending only on that parameter, its likelihood can be formed and maximized to produce a second stage MLE
(Savage, 1976, p.455, footnote 20). 
 Fisher (1915, 1921) used the procedure
 to estimate the  variance 
and the correlation coefficient
of  normal population 
 but  has  never formulated this `` second criterion''.
He  has  never discussed  the relationship between the original and the
second criterion, why he preferred the latter in  1912-1921 and changed his mind in 1922
(Aldrich, 1997,  p. 166, left column, lines 22-35).
There were  neither motivating theory nor details for the implementation of the 2-stage procedure. For example,
which estimate
to choose if the second step estimate has smaller bias but larger MSE than the estimate obtained in the first step? Why is the estimate in the second step better than that in the first  step?  

MLP was  introduced
 by 
Fisher $(1922, 1925)$ who established  asymptotic optimality of  the
MLE  $\hat \theta$
of  $\theta$  for
various $x$-models. 
The notions of the first and second order efficiency of an estimate
revealed asymptotic optimality properties of 
$\hat \theta$ 
(Rao, 1962, Efron, 1975).
A decision theoretic approach
showed that $\hat \theta$
is finite sample efficient with respect to the 
mean squared error
 of the scores and within a large class 
of estimates  (Yatracos, 1998).

 Nevertheless,
several examples
 in the literature showed  that the MLE
is either biased, or  inconsistent, or there are 
better estimates. 
Many of the examples and  criticisms 
appear  in LeCam (1990) who added ``It might simply mean we have not yet translated into mathematics the basic principles 
which underlined Fisher's intuition.''  A lot of research was devoted   to relax the  criticisms by providing  MLE's corrections 
thus violating  MLP that did not advocate correction.
 Firth (1993) observed that 
most  methods are corrective in character rather than preventive, i.e. the MLE   is first calculated and then 
corrected,  and proposed a preventive approach with  systematic correction of the likelihood  equations (LEs).


 This work is
 motivated  from several  
 MLEs for the shape parameter $\psi$  that are unbiased only when the location
$\theta$ is known.
The goals are:\\ $a)$ to examine whether there is a theoretical explanation for this phenomenon, $b)$  to  correct the bias adhering to MLP. 

The obtained results for  $a)$
 show that $\hat \psi$'s bias in these examples is not a coincidence and indicate how to achieve $b)$
 by not adhering to  Fisher's model  specification 
approach (Fisher, 1922, 1925)  that dictates to 
determine {\em once and  for all} from the data  the population model used to obtain the LEs. 
Fisher's approach  indirectly  implies that  the stochastic quantities in  the LEs 
have the same information with the data.
However,
when  $\theta$ is replaced by $\hat \theta$ 
 in the LE to be solved for $\psi$ a new situation arises. This modified  LE 
has a new stochastic component and  the 
{\em updated data} $Y$ in it  introduces   inaccuracy with respect to the original  LE because {\em i)}  $\theta$ 
is replaced by $\hat \theta$ and {\em ii)} $Y$'s degrees of freedom change. For example, with a sample $x=\{X_1,\ldots, X_n\}$  from the normal model with mean $\theta$ and variance
 $\sigma^2$ the LE for $\sigma^2$ depends on $\sum_{i=1}^n (X_i-\theta)^2$ that has $n$ degrees of freedom.
 When the MLE
$\bar X$
replaces 
$\theta$ inaccuracy is introduced and the ``updated data'', $Y,$    in the  LE used to obtain 
$\hat \sigma^2$  is 
$$Y=\sum_{i=1}^n (X_i-\bar X)^2.$$ 
This new LE is not that of  a 
$\chi^2$-distribution with $n-1$ degrees of freedom, i.e. $Y$'s distribution, thus it ``does not correspond to a proper model''.


The  proposed preventive approach
suggests to replace the 
LE to be solved for  $\psi$ 
after plugging $\hat \theta$ in it  with the LE from $Y$'s distribution,
 thus adhering to MLP. The data $Y$ is a multiple of MLE 
$\hat \psi$ used in the 2-stage procedure.
Using  {\em model  updated}  LEs
 unbiased 
$\hat \psi_{MU}$ are
obtained for the shape parameters of the  normal and the shift-exponential models; the variance estimate $\hat \psi_{MU}$
for the Neyman and 
Scott (1948)
problem is unbiased and consistent;  the shape parameter's estimate $\hat \psi_{MU}$  for the Pareto distribution improves both the bias and the variance of $\hat \psi$  and, in addition,  by  parametrizing the model  with $\psi^{-1}$ 
its MUMLE  is unbiased 
contrary to the  MLE.

 MUMLE's approach
justifies from a frequentist's view the likelihood correction in the
MML estimation method (Wallace and Boulton, 1968, Wallace and Freeman, 1987, Wallace, 2005)
and in the Minimum Description Length Criterion (Rissanen 1984, 1987). 
Both methods
assume a prior distribution 
 but have different philosophy for its choice and use (Rissanen, 1987, p. 226, Wallace and Freeman, 1987, p. 251).
Model update satisfies one of
Rissanen's criticisms for the MLE ``... the maximized likelihood $P(x|\hat \theta(x))$ no longer defines a proper distribution'' (1987, p. 224).

MUMLE's  formulation violates Fisher's model  specification approach  but  adheres to  MLP and
more precisely to MUMLP.  MUMLE should be  explored further.
The  2-stage procedure
does not adhere to MLP which does not allow for corrections.
 It is  a bias corrective approach 
that does not touch the heart of the matter, i.e.,  
it does not explain why the difference in bias 
 occurs
and does not  motivate the remedy.
These are the reasons we prefer the formulation for the MUMLE approach.
The  puzzling question is  Fisher's  rigidity with  the  model specification. A possible explanation  is the  Bayesian flavor involved with model updating.
 
\section {MLE's  Bias Pathology}

 \quad Let  the data $x$ be a  random vector in $R^d$ 
having density $f(x|\theta, \psi)$ with respect to Lebesgue measure,
   parameters $ \theta \in R^p, \ \psi  \in R$ all unknown and  
with  the $\psi$-score
 $U_{\psi}$ satisfying
\begin{equation}
\label{eq:psiscore}
U_{\psi}(x,\theta, \psi)=\frac{\partial{\log f(x|\theta, \psi)}}{\partial {\psi}} \neq 0 \  a.s. \ f(\cdot|\theta,\psi),
\end{equation}
 \begin{equation}
\label{eq:unique}
\forall  \ x, \theta, \hspace{2ex}  U_{\psi}(x,\theta, \psi)=0 \hspace{2ex}  \mbox{has unique solution,}
\end{equation}
\begin{equation}
\label{eq:modelreg}
 E_{\theta, \psi}U_{\psi}(x, \theta, \psi)=0 \ (\psi \mbox{-regularity});
\end{equation}
$E_{\theta, \psi}$ denotes expectation with respect to $f(x|\theta, \psi), \ d \ge 1, \ p\ge 1. $

Assume that  MLE  $\hat \theta$ of $\theta$ and $U_{\psi}$ are  used to obtain    MLE  $\hat \psi$ 
 such that 
\begin{equation}
\label{eq:1}
U_{\psi}(x, \hat \theta, \hat \psi)=0.
\end{equation}
 

It is seen 
in Proposition \ref{p:nsc} $a)$
that $\psi$-regularity  (\ref{eq:modelreg}) may most often cause bias for   $\hat \psi$  because it 
is expected 
to imply
that
$E_{\theta, \psi}U_{\psi}(x, \hat  \theta, \psi)$ does not vanish, especially if 
$\theta$'s dimension $p$ is large. Using instead  the score
for the data $Y$ (i.e. $\hat \psi$) to determine $\hat \psi_{MU}$  
 this
drawback
is avoided for some models thus motivating  the use of MUMLE.


\bep \footnote{Proofs are in the Appendix.}
 (MLE's inherent bias pathology)
\label{p:nsc}
 Let $x$ be data in $R^d$  from 
$f(x|\theta, \psi)$  with $\theta \in R^p, \ \psi \in R$ both unknown
with  the $\psi$-score
$ U_{\psi}$
 satisfying (\ref{eq:psiscore})-(\ref{eq:modelreg}) and 
$\hat \psi$  obtained from (\ref{eq:1}); $\hat \theta$ is the MLE of $\theta, \ d \ge 1, \   p\ge1.$

\quad  a) If  
 $ \frac { \partial {U_{\psi}(x, \hat \theta,\psi)} }{ \partial {\psi}}=C$ 
  is  fixed constant,
$ C \neq 0,$
  $\hat \psi $ is biased estimate of $\psi$ if and only if
\begin{equation}
\label{eq:nsc1}
 E_{\theta, \psi}U_{\psi}(x, \hat \theta, \psi) \neq  0
\end{equation} 
at least for $\psi=\psi_0.$
Since  (\ref{eq:modelreg}) holds $\hat \psi$ is expected to be biased.

 b) If  $ \frac { \partial {U_{\psi}(x, \hat \theta,\psi)} }{ \partial {\psi}}=C(x,\hat \theta, \psi) $  exists in a neighborhood of 
$\psi_0,$
$\hat \psi $ is biased estimate of $\psi$ if and only if
\begin{equation}
\label{eq:nsc2}
 E_{\theta, \psi}\frac{U_{\psi}(x, \hat \theta, \psi)}{C(x,\hat \theta, \psi^*)} \neq  0
\end{equation} 
at least for $\psi=\psi_0; \ \psi^*$ is  between $\hat \psi$ and $\psi_0.$
$\hat \psi$ is expected
 to be biased.



\enp

A simple result follows  motivating the use of MUMLE  when $Y$'s distribution depends only on $\psi.$
 
\bec
Under the assumptions of  Proposition \ref{p:nsc} $a)$  but with $\psi$ the only model parameter,
$\psi$-regularity (\ref{eq:modelreg})  implies  that  
$\hat \psi$ is unbiased for $\psi.$
\enc

The next proposition can be  used to show $\hat \psi$ is biased.

\bep
\label{c:tool}
Let $T(x,\theta, \psi)$ be a functional for which 
(\ref{eq:1}) holds with $T$ instead of $U_{\psi}$,
$ \frac { \partial {T} } { \partial {\psi}}$  is a constant $C (\neq 0)$ 
and for
 $\psi_0$ it holds
$$E_{\theta, \psi_0}T(x,\hat \theta, \psi_0) \neq 0.$$
Then $\hat \psi$ is biased estimate of $\psi.$
\enp

 
 
When $U_{\psi}(x, \theta, \psi)$ has the form
\begin{equation}
\label{eq:abuse}
U_{\psi}(x, \theta, \psi)=\frac{ U^* (x, \theta, \psi)}{\tilde h(\psi)},
\end{equation}
 (\ref{eq:unique})-(\ref{eq:1}) hold also for  $U^*;$
$\tilde  h$ is a real valued function. The equation to be solved for $\psi$ has the form
\begin{equation}
\label{eq:formustar}
U^*(x,\theta, \psi)=C(x,\theta)\psi +D(x,\theta)=0.
\end{equation}
 $U^*$  is a useful tool
that  will play the role of $T$ when applying Proposition \ref{c:tool}.\\

With the next proposition $\hat \psi$'s bias is confirmed {\em directly} for some models.
 
\bep
\label{p:elab1}
For $f(x|\theta, \psi)$ with $\hat \theta$ the MLE for $\theta$ assume in addition to (\ref{eq:psiscore})-(\ref{eq:modelreg})  that \\
a)  $\psi>0,$\\
b) 
 \begin{equation}
\label{eq:specf}
\log f(x|\theta, \psi)=\frac{C}{A} \log \psi - \frac{D(x, \theta)}{A \psi}+ g(x)
\end{equation}
which implies  that 
 \begin{equation}
\label{eq:specU}
U_{\psi}(x|\theta, \psi)=\frac{C\psi+D(x,\theta)}
{A\psi ^2};
\end{equation}
$C$ is a constant, $D$ is a  function with positive values, $A>0$  and $g$ is a real valued function of $x.$\\
 Then, $\hat \psi$ is biased for $\psi.$  
\enp



Proposition \ref{p:elab1} is used in  Examples 2.1-2.4.

\beex {\em  Let $x=\{X_1,\ldots, X_n\}$ be
  i.i.d. normal  random variables with mean $\theta$ and variance $\psi.$ Then
$f(x|\theta, \psi)$ satisfies (\ref{eq:specf}), $\hat \theta=\bar X$  and  $U_{\psi}$ has  form (\ref{eq:specU}) with 
$$C=-n, \  D(x,\theta)=\sum_{i=1}^n (X_i-\theta)^2, \ A=2.$$
From Proposition \ref{p:elab1}  $\hat \psi$ is biased for  $\psi.$}
\enex

\beex (The Neyman-Scott problem) {\em Let $\{X_{ij},  j=1,...,k \}$ be a sample from a normal distribution with 
mean $\theta_i$ and 
variance
 $\psi , i=1,...,n,$ and let $x$ represent all the observations.  The samples are independent and $\hat \theta_i=\bar {X}_i,$ 
$i=1,\ldots, n.$ 
Then $f(x|\theta, \psi)$ satisfies (\ref{eq:specf}) 
and $U_{\psi}$  has form (\ref{eq:specU}) with 
$$C=-nk, \  D(x,\theta_1,\ldots, \theta_n)=\sum_{i=1}^n \sum_{j=1}^k (X_{ij}-\theta_i)^2, \ A=2.$$
From Proposition \ref{p:elab1} it follows that  $\hat \psi$ is biased for  $\psi.$}
\enex

\beex {\em Let $x=\{X_1,\ldots, X_n\}$ be 
 i.i.d.  random variables from a shifted exponential density  $f$  with parameters 
$\theta$ and $\psi \ (>0),$ 
\begin{equation}
\label{eq:expo}
f(w,\theta, \psi)=\psi^{-1}e^{-(w-\theta)/\psi}I_{[\theta, \infty)}(w);
\end{equation}
$I$ denotes the indicator  function.
Then $f(x|\theta, \psi)$ satisfies (\ref{eq:specf}), $\hat \theta$ is the smallest observation $X_{(1)}$ and $U_{\psi}$  has
 form (\ref{eq:specU}) with 
$$C=-n, \  D(x,\theta)=\sum_{i=1}^n (X_i-\theta), \ A=1.$$
From Proposition \ref{p:elab1}   $\hat \psi$ is biased for  $\psi.$}
\enex

\beex \label{ex:pareto} (Pareto family with  non-usual parametrization of  the shape parameter.) {\em Let $x=\{X_1,\ldots, X_n\}$ be 
  i.i.d.  random variables with density
\begin{equation}
\label{eq:paretorepar}
f(w|\theta, \psi^*)=\frac{1}{\psi^*} \theta^{1/\psi^*}w^{-(\frac{1}{\psi^*}+1)}I_{[\theta, \infty)}(w), \ \psi^*>0, \ \theta >0;
\end{equation}
$I$ denotes the indicator function. Then $f(x|\theta, \psi)$ satisfies (\ref{eq:specf}), 
$\hat \theta$ is the smallest observation $X_{(1)}$ and $U_{\psi^*}$  has
form (\ref{eq:specU}) with 
$$C=-n, \  D(x,\theta)=\sum_{i=1}^n \log  \frac{X_i}{\theta}, \ A=1. $$
From Proposition \ref{p:elab1}  MLE $\hat \psi^*$ is biased for  $\psi^*.$}
\enex

The proposition that follows presents conditions under which $\hat \psi$  is biased.
The definition of a complete family of densities is provided according to  Lehmann and Scheff\'e (1950). 

\bef
Let ${\cal G}=\{g(u|\eta), \ \eta \in {\cal H}\}$ be a family of densities of a random variable (or statistic) $U$  indexed by 
the parameter set ${\cal H}.$  ${\cal G}$  is complete if for any function $\phi$ satisfying 
$$E_{\eta}\phi(U)=0 \ \forall \ \eta \in {\cal H}$$
it holds that $\phi(u)=0$ for every $u$ except for a set of $u$'s having probability zero for all $\eta \in {\cal H}.$
\enf

 \bep 
\label{p:nsc2}
a) Under the assumptions and the notation  of Proposition \ref{p:nsc} a), with  $C(x,\hat \theta, \psi)$ a constant $C$ and
\begin{equation}
\label{eq:ascondition}
f(x|\theta, \psi)>0 \ \forall \  x \in U \subset R^d, \forall \ \theta, \psi,
\end{equation}
 if  the family $\{ f(x, \theta, \psi), \theta \in R  \}$ is complete for each fixed $\psi$  and the distribution of $U_{\psi}(x|\hat \theta, \psi)$ depends also on $\theta,$   then 
$\hat \psi$ is  biased estimate of $\psi.$

b)  Under the assumptions  and the notation of  Proposition \ref{p:nsc} b),  for  general  $C(x,\hat \theta, \psi)$ existing
in neighborhoods of $\psi_0$ and $\tilde \psi_0$  and
with (\ref{eq:ascondition}) holding,
 if
 the family $\{ f(x, \theta, \psi), \theta \in R \}$ is complete for each fixed $\psi$ and
the distribution of $\frac{U_{\psi}(x|\hat \theta, \psi)}{C(x,\hat \theta, \psi*)}$  for $\psi=\psi_0, \tilde \psi_0,$ depends also on $\theta,$  then  
 $\hat \psi$ is biased.
\enp


\ber {\em Proposition  \ref{p:nsc2} motivates the use of MUMLE and the 2-stage procedure when $\hat \psi$'s distribution does not
depend on $\theta.$
 Proposition \ref{p:nsc2} a) does not apply in Examples 2.1-2.4  because
$U_{\psi}(x|\hat \theta, \psi_0)$'s  distribution
does not depend on $\theta.$}
\enr

\section{Fisher's specification problem, MUMLE and the  MML method}

\quad According to Fisher(1922): ``...
 The data is to be replaced by few quantities that will 
contain as much 
as possible of the relevant information contained in the original data. 
This object is accomplished 
by constructing   a hypothetical infinite population  of which the actual 
data are regarded as constituting 
a random sample{\em (the specification problem).} ...
 The problems of specification are entirely a matter 
for the practical statistician. The discussions of theoretical statistics 
may
be regarded as alternating between problems of estimation and problems of 
distribution.''

 We include the specification problem in these alternating discussions.
The goal is that the $k$-th LE to be solved, $k \ge 2,$ maximizes a proper likelihood, i.e. a likelihood that coincides with that 
of the data $Y$  in it after replacement of other parameter values with their MLEs.
Results in section 2 suggest that bias may be reduced.\\
{\em The MUMLP approach:} 
Let $f(x|\theta_1,...,\theta_p)$ 
be  the  density of the data $x; \theta_1,...,\theta_p$ 
are real valued parameters.  
Assume that $k-1$ likelihood equations have been solved obtaining 
 estimates $\hat \theta_1,..., \hat \theta_{k-1},$ 
respectively,  of $\theta_1,...,\theta_{k-1}, k -1< p.$
The LE for $\theta_k$ has form (\ref{eq:formustar}) with $\theta_k$ instead of $\psi,$  and solving 
it we obtain 
$$\hat \theta_k=-\frac{D(x,\hat \theta_1,\ldots, \hat \theta_{k-1})}{C(x,\hat \theta_1,\ldots, \hat \theta_{k-1})}=Y.$$
When $Y$'s density 
depends only on $\theta_k$ it 
is used as 
model to obtain MUMLE $\hat \theta_{k,MU}.$\\

In the examples presented in the next section 
the distribution of $Y$ is easy to obtain. 
If 
$Y$'s  distribution  is not immediately accessible, as in the case of a sample $x=\{X_1,\ldots, X_n\}$  from a Gamma distribution with  two unknown parameters and $Y=\Pi_{i=1}^n X_i /\bar X_n^n,$ other methods can be used to obtain a
LE  from a proper model. One possibility is to use the machinery of the MML87 method (Wallace and Freeman, 1987, Wallace, 2005) 
for the model $f(x|\theta)$ with prior $h(\theta)$
and choose, according to a criterion, one of the estimates obtained from a data-dependent class of priors.
\\

{\em The MML87 method: } The MML estimate of $\theta (\in R^p)$  is the value $\hat {\theta}_{MML}$ maximizing
\begin{equation}
\label{eq:mml871}
\log h(\theta) +\log f(x|\theta)-\frac{1}{2} \log |I_x (\theta)|;
\end{equation}
$h(\theta)$ is a prior and $|I_x(\theta)|$ is  the determinant of the Fisher's information matrix  for  $x, \  p\ge 1.$\\


The next propositions motivate  the use 
of the MML approach for frequentist inference.

\bep
\label{p:mmlfirth}
If  $\theta (\in R^p)$ are the canonical parameters of an exponential family model, the MML estimates remove the $O(n^{-1})$ term in $\hat {\theta}$'s  bias when 
\begin{equation}
\label{eq:firthcanon}
h(\theta) \propto |I_x(\theta)|.
\end{equation}
\enp

\ber {\em Proposition \ref{p:mmlfirth} can be extended for exponential family models in non-canonical parametrization as well as for non-exponential models with the proper choice of $h(\theta)$ along the lines in  Firth (1993, p. 30, sec. 4).} 
\enr

The  proposition that follows provides conditions for a model with parameters $\theta$ and $\psi$ and $\hat \psi$ function of 
$\hat \theta$ under which  the MUMLE estimates  $\hat \theta, \ \hat \psi_{MU}$ coincide with MML estimates  $\hat \theta_{MML}, \hat \psi_{MML}.$

\bep
\label{p:mlemcmml}
Assume that  the data $x$ has density $f(x|\theta, \psi), \ \theta \in R^p, \ \psi \in R,$   that MLEs  $\hat \theta, \ \hat \psi$ are obtained,
$\hat \psi$ is a function of $\hat \theta$  and  $Y(i.e. \  \hat \psi)$ has density $g_Y(y|\psi).$
Assume in addition that \\
a) $|I_x(\theta, \psi)|=|I_x(\psi)|,$ \\
b) there are functions $\phi(\psi), \ u(y):$
\begin{equation}
\label{eq:likg}
\log f(x|\hat \theta, \psi)-\log g_Y(y|\psi) =  \log \phi(\psi) +u(y).
\end{equation}
Then, MML estimates $\hat \theta_{MML}$ and $\hat \psi_{MML}$  coincide, respectively, with $\hat \theta$ and
$\hat \psi_{MU}$  if the prior
\begin{equation}
\label{eq:mcprior}
h(\theta,\psi) \propto \frac{|I_x(\psi)|^{1/2}}{\phi(\psi)}.
\end{equation}
\enp

\ber {\em The assumptions in Proposition \ref{p:mlemcmml} hold at least under the set-up of Example 2.1 for which
$$\phi(\psi) \propto \psi^{-1/2}, \hspace{5ex}   |I(\theta, \psi)|=|I(\psi)|=2n^2/\psi^2.$$
Then,
$$h(\theta, \psi) \propto \psi^{-1/2}$$
that is the prior used to obtain $\hat \theta_{MML}, \ \hat \psi_{MML}$ (Wallace, 2005, p. 250).}
\enr

 \section{Examples-MUMLE's Applications}

An elementary Lemma follows  to be used in the examples.

\bel  
\label{l:thelemma}  Let $W$ be a chi-square random variable with $k$ 
degrees
of freedom and let $Y= W \tau^2, \tau>0$. Then,\\ a) 
$Y$'s density
has the
form $C_k \exp\{-y/2\tau^2\} y^{(k-2)/2} \tau^{-k}, C_k (>0)$ is a 
constant.\\
b) The likelihood equation, corresponding to    Y is 
$$-k \tau^2 + Y = 0$$  
 and
the MLE $\hat {\tau}^2$ is given by $Y/k.$
\enl

 The first example is the variance estimation problem for a normal sample with unknown mean.
The MUMLE of the variance is its unbiased estimate that is also the MML estimate (Wallace and Boulton, 1968, p.190) and
Firth's (1993, p. 34, l. 1) bias corrected estimate.

\beex (Example 2.1 continued) {\em The LE for $\psi$
with $\hat \theta=\bar X$ is  
$$ -n \psi + \sum_{i=1}^n 
(X_i-\overline {X})^2 = 0, \hspace{5ex} Y = \sum_{i=1}^n (X_i-\overline {X})^2$$
 and $Y$'s distribution  follows
 from  Lemma \ref{l:thelemma}  with $\tau$ and $k$ taking,respectively, values $\sqrt{\psi}$ and $n-1.$ 
The model updated LE is
 $$-(n-1) \psi + \sum_{i=1}^n (X_i-\overline{X})^2 = 0.$$
 The $MUMLE$ of $\psi$  is its $UMVU$
estimate
 $$ (n-1)^{-1} \sum_{i=1}^n (X_i-\overline{X})^2.$$}
\enex

\beex (Example 2.2 continued, the Neyman-Scott problem) {\em 
 The LE for $\psi$ after replacing $\theta_i$ by its MLE $\bar X_i$ (for every $i$)  is 
$$-nm\psi +\sum_{i=1}^n\sum_{j=1}^m(X_{ij}-\bar {X}_i)^2 = 0, \hspace{5ex} 
Y = \sum_{i=1}^n\sum_{j=1}^m(X_{ij}-\bar {X}_i)^2. $$  
Using  $Y$'s model from Lemma \ref{l:thelemma}  with $ k=n(m-1),$  the MUMLE  is 
$$n^{-1}(m-1)^{-
1}\sum_{i=1}^n\sum_{j=1}^m(X_{ij}-\bar {X}_i)^2,$$ an unbiased 
and consistent estimate of $\psi.$}
\enex

For the Neyman-Scott problem one of Firth's (1993, p. 35) estimates of $\sigma^2,$  $A^{(O)},$ is  unbiased and consistent  while
the other estimate, $A^{(E)},$ is consistent.  The  MML estimate  obtained is  consistent and asymptotically unbiased (Dowe and Wallace, 1997,  p. 617, Wallace, 2005, p. 202).

\beex  (Example 2.3 continued) {\em 
$\hat \theta$ is the smallest observation   $X_{(1)}$ and
the LE  for $\psi$ 
 is 
$$-n \psi+ \sum_{i=1}^n(X_{(i)}-X_{(1)})=0, \hspace{5ex} Y=\sum_{i=1}^n(X_{(i)}-X_{(1)}).$$
$Y$ follows  Gamma distribution with parameters $\psi$ and
$n-1.$ The LE  for $Y$ is
$$-(n-1) \psi+ \sum_{i=1}^n(X_{(i)}-X_{(1)})=0$$
and the MUMLE of $\psi$ is
$$\frac{\sum_{i=1}^n(X_{(i)}-X_{(1)})}{n-1}$$ 
that is also the UMVU estimate  (Arnold, 1970, p. 1261).}

\enex
  
In the Pareto family example that follows with  parameters $\psi$ and $\theta$ both unknown
$\hat \psi_{MU}$ 
reduces  by 50\% the bias  of the MLE $\hat \psi$  and has also smaller variance. With this parametrization $\hat \psi$ is not unbiased
even when $\theta$ is known.  Using the parametrization $\psi=1/\psi^*,$  MLE  $ \hat \psi^*$  is unbiased for $\psi^*$ when $\theta$
is known but when $\theta$ is unknown MUMLE  $\hat \psi^*_{MU}$ is unbiased. 
\beex {\em 
 Let 
$X_1, \cdots, X_n$
be independent random variables from 
Pareto 
 density (\ref{eq:paretorepar}) with $\psi^*=\psi^{-1}, \ \psi>0.$
The log-likelihood 
of the sample is
$$
n \log \psi + n \psi \log \theta -(\psi+1) \sum_{i=1}^n \log X_i  + \sum_{i=1}^n \log I_{[\theta, \infty)}(X_i)
$$
and  $\hat \theta$ is the smallest observation, $X_{(1)}.$
The score and the MLE are, respectively,
$$U_{\psi}(X, \hat \theta, \psi)
= n- \psi \sum_{i=2}^n \log \frac{X_i}{X_{(1)}}, \hspace{5ex}
\hat \psi=\frac{n}{\sum_{i=2}^n \log \frac{X_i}{X_{(1)}}}.$$ 
Since 
$$Y=\sum_{i=2}^n \log \frac{X_i}{X_{(1)}}$$
 has a $\Gamma(n-1, \psi)$ distribution (see, e.g,  Baxter, 1980,
 p. 136, l. -6  and references therein)
$\hat \psi$ is biased and 
$$ \ E\hat \psi-\psi=\frac{2}{n-2}\psi, \ 
Var(\hat \psi)=\frac{n^2}{(n-2)^2(n-3)}\psi^2.$$

The  updated  score based on the data $Y$ and  MUMLE $\hat \psi_{MU}$ are, respectively,
$$(n-1)-\psi Y, \hspace{5ex}
\hat \psi_{MU}=\frac{n-1}{\sum_{i=2}^n \log \frac{X_i}{X_{(1)}}},$$
with
$$ E\hat \psi_{MU}-\psi=\frac{1}{n-2}\psi,
 \   Var(\hat \psi_{MU})=\frac{(n-1)^2}{(n-2)^2(n-3)}\psi^2. $$
Observe that $\hat \psi_{MU}$ improves both the bias and the variance of $\hat \psi.$\\

Using instead density  (\ref{eq:paretorepar}) the $\psi^*$-score and the MLE are, respectively,
$$U_{\psi^*}(X, \hat \theta, \psi^*)
=- n \psi^*+\sum_{i=2}^n \log \frac{X_i}{X_{(1)}}, \hspace{5ex}
\hat \psi^*=\frac{\sum_{i=2}^n \log \frac{X_i}{X_{(1)}}}{n}.$$
$\hat \psi^*$ is biased; see Example \ref{ex:pareto}.
Using the model from data $$Y=\sum_{i=2}^n \log \frac{X_i}{X_{(1)}}$$  the updated  score 
and  MUMLE $\hat \psi^*_{MU}$  are, respectively
$$-(n-1)\psi^*+Y, \hspace{5ex}
\hat \psi^*_{MU}=\frac{\sum_{i=2}^n \log \frac{X_i}{X_{(1)}}}{n-1}.$$
$\hat \psi^*_{MU}$  is unbiased for  $\psi^*.$ 


}
\enex

  m hospitality during my summer visits when most of the results were obtained. This research was partially supported CUT.

\begin{center}
{\bf Appendix} 
\end{center}

 {\bf Proof  of Proposition \ref{p:nsc}:} $ a)$ 
Make a  Taylor expansion of $U_{\psi}(x, \hat \theta,\hat  \psi)$ around $ \psi$ using $U_{\psi}$'s linearity in $\psi,$ 
\begin{equation}
\label{eq:exp}
U_{\psi}(x, \hat \theta, \hat  \psi)=U_{\psi}(x, \hat \theta,  \psi)+(\hat \psi- \psi)C.
\end{equation}

From (\ref{eq:1}) it follows that 
$$E_{\theta,\psi}(\hat \psi- \psi)=-C^{-1} E_{\theta,\psi} U_{\psi}(x, \hat  \theta, \psi)  \neq 0$$
 if and only if $E_{\theta,\psi} U_{\psi}(x, \hat  \theta, \psi)\neq 0. $ 

$ b)$  Equation (\ref{eq:exp}) remains valid with $C=C(x,\hat \theta, \psi)$ 
evaluated at $\psi=\psi^*$ between $\psi$ and $\hat \psi.$
 Then $\hat \psi$ is biased if and only if
\begin{equation}
\label{eq:expran}
E_{\theta,\psi} U_{\psi}(x, \hat \theta, \psi) C^{-1}(x,\hat \theta, \psi^*) \neq 0.
\end{equation}
Most often (\ref{eq:expran}) will hold. To examine this expectation further 
make   a second order Taylor approximation of the left  side in (\ref{eq:expran}) around $E_{\theta, \psi}U_{\psi}(x, \hat \theta, \psi)$ (denoted by $EU_{\psi}$) and $E_{\theta,\psi}C(x,\hat \theta, \psi^*)$ (denoted by $EC$) assuming negligibility of the
remainder, 
\begin{equation}
\label{eq:Texp}
E_{\theta,\psi}\frac{U_{\psi}}{C} \approx \frac{EU_{\psi}}{EC} -\frac{Cov(U_{\psi}, C)}{E^2C}+\frac{Var(C)EU_{\psi}}{E^3C}.
\end{equation}
Whether or not $EU_{\psi}=0,$ (\ref{eq:Texp}) is not expected to vanish. $\hspace{3ex} \hfill \Box$

{\bf Proof  of Proposition \ref{c:tool}:} Follows along the proof of Proposition \ref{p:nsc}{\em  a)}  with $T$ instead of $U_{\psi}$
since $T$ is linear in $\psi.$ $\hspace{3ex} \hfill \Box$

{\bf Proof of Proposition \ref{p:elab1}:}  From   (\ref{eq:specf}) 
 \begin{equation}
\label{eq:equivDf}
f(x|\hat \theta, \psi) >  f(x|\theta, \psi) \ \forall  \ \psi \Leftrightarrow D(x,\hat \theta) < D(x,\theta).
\end{equation}
Thus,  from (\ref{eq:specU}) it follows that
$$   U_{\psi}(x|\hat \theta, \psi) < U_{\psi}(x| \theta, \psi) \  a.s. $$
$$\Rightarrow \  E_{\theta, \psi}U_{\psi}(x|\hat \theta, \psi)< E_{\theta, \psi}U_{\psi}(x| \theta, \psi)=0$$
 from (\ref{eq:modelreg}).
From  (\ref{eq:specU}) it also holds that 
$$E_{\theta, \psi}[C\psi+D(x,\hat \theta)] \neq 0$$
and from  Proposition \ref{c:tool}  with
$$T(x, \theta, \psi)=C\psi+D(x,\theta)$$
$\hat \psi$ is biased. $ \hspace{3ex} \hfill \Box$

{\bf Proof of Proposition \ref{p:nsc2}:}
$a)$  The result is proved by contradiction. Assume that $\hat \psi$ is unbiased.
Then from Proposition \ref{p:nsc} $a)$   for $\psi_0$ 
\begin{equation}
\label{eq:nsccontradiction}
E_{\theta, \psi_0} U_{\psi}(x, \hat \theta, \psi_0)=0  \ \forall \ \theta.
\end{equation}
Let 
$$K(\psi_0)=\{x:U_{\psi}(x, \hat \theta, \psi_0)=0\}.$$
Since $ U_{\psi}(x, \hat \theta, \psi_0)$ is  function of $x$  only, by assumption  its distribution
depends on both $\theta$ and $\psi_0$ and  the family $\{f(x|\theta, \psi_0), \theta \in R\}$ is complete,
 it follows from (\ref{eq:nsccontradiction}) that 
\begin{equation}
\label{eq:contrac2}
P_{\theta, \psi_0}[U_{\psi}(x, \hat \theta, \psi_0)=0]=P_{\theta, \psi_0}[K(\psi_0)]=1  \ \forall \ \theta.
\end{equation} 
Equalities (\ref{eq:contrac2}) hold also for 
 $\tilde \psi_0 \neq \psi_0$ and for $x\in K(\psi_0) \cap K(\tilde \psi_0) (\neq \emptyset) $ 
the likelihood equation for $\psi$  has 2 solutions,  $\psi_0$ and  $\tilde \psi_0,$
leading to contradiction because of (\ref{eq:unique}).


$b)$  Assume that $\hat \psi$ is unbiased.
From
Proposition \ref{p:nsc} $b)$   for $\psi_0$ it holds
$$E_{\theta, \psi_0} \frac{U_{\psi}(x, \hat \theta, \psi_0)}{C(x, \hat \theta, \psi^*)}=0 \   \forall \ \theta.$$
Since $  \frac{U_{\psi}(x, \hat \theta, \psi_0)}{C(x, \hat \theta, \psi^*)}$
 is function of $x$ only, its distribution depends on both $\theta$ and $\psi_0$  and family $\{f(x|\theta, \psi_0), \theta \in R \}$ is complete it follows that 
$$  \frac{U_{\psi}(x, \hat \theta, \psi_0)}{C(x, \hat \theta, \psi^*)}=0 \ a.s.$$
which implies that 
 $$P_{\theta, \psi_0}[U_{\psi}(x, \hat \theta, \psi_0)=0]=1 \ \forall \ \theta.$$
The proof follows as in part {\em a)}. $ \hspace{3ex} \hfill \Box$


{\bf Proof of Proposition \ref{p:mmlfirth}:} Replacing (\ref{eq:firthcanon}) in (\ref{eq:mml871}) it follows that $\hat {\theta}_{MML}$
is the value maximizing 
$$\log f(x|\theta)+\frac{1}{2} \log |I_x (\theta)|$$
and the result follows from Firth (1993, p. 30, sec. 3).  $ \hspace{3ex} \hfill \Box$

{\bf Proof of Proposition \ref{p:mlemcmml}:}
Replacing $h$ from (\ref{eq:mcprior})  in (\ref{eq:mml871}) the MML  log-likelihood is
\begin{equation}
\label{eq:likmlemcmml}
c-\log \phi(\psi)+\log f(x,\theta, \psi);
\end{equation}
$c$ is a constant. It follows that 
$$\hat \theta_{MML}=\hat \theta.$$
From (\ref {eq:likg}) and (\ref{eq:likmlemcmml}) the  MML log-likelihood for $\psi$ is 
$$c-\log \phi(\psi)+\log f(x,\hat \theta, \psi)=c+ \log g_Y(y|\psi) + u(y)$$
and 
$$\hspace{32ex} \hat \psi_{MML}=\hat \psi_{MU}. \hspace{25ex} \hfill \Box$$

{\bf Proof of Lemma \ref{l:thelemma}: } The density of $W$ is given by $C_k w^{(k-2)/2} \exp\{-
w/2\},
C_k>0.$
Thus, the density of $Y$ is $C_k \exp\{-
y/2\tau^2\}
 (y/ \tau^2)^{(k-2)/2} \tau^{-2}. \hspace{3ex} \hfill \Box $

\end{document}